\documentclass[12pt]{article}

\usepackage{amsfonts,amssymb,eucal}
\begin{document}

\title{{\large \bf THE MODULI SPACE OF HYPERBOLIC}\\{\large \bf CONE 
STRUCTURES}}

\author{{\sc Qing Zhou}\thanks{The project is partially supported by 
grants of NNSFC, FEYUT of SEDC and HYTEF.} \\
{\small{\sc Department of Mathematics, East China Normal University}}\\
{\small{\sc Shanghai, 200062, China}}}

\maketitle

\vspace{0.1in}

{\bf Introduction.} Roughly speaking, a cone structure is a 
manifold together with a link each of whose component has a 
cone angle attached. It is a kind of singular manifold 
structure. If each cone angle is of the form $2\pi /n$, 
for some integer $n$, the cone structure becomes an orbifold
structure. Unlike an orbifold structure, the cone structure is 
not a natural concept, but it turns out to be very important 
in the study of geometric orbifolds.

In this paper, we will consider 3-dimensional geometric cone 
structures. The main result in this paper is an existence and 
uniqueness theorem for 3-dimensional hyperbolic cone structures.

{\sc Theorem A.} {\em Let $\Sigma$ be a hyperbolic link with 
$m$ components in a 3-dimensional manifold $X$. Then the 
moduli space of marked hyperbolic cone structures on the pair 
$(X,\Sigma)$ with all cone angles less than $2\pi /3$ is 
an $m$-dimensional open cube, parameterized naturally
by the $m$ cone angles.}

Theorem A is an analogue of Mostow's rigidity theorem. If we 
have two hyperbolic cone structures $C_1$ and $C_2$ with cone
angles less than $2\pi /3$, then $C_1$ and $C_2$ are isometric 
if and only if there is a homeomorphism between $(X_1,\Sigma _1)$ 
and $(X_2,\Sigma _2)$ so that corresponding cone angles are the same.
The proof of the theorem goes in a similar way as Thurston 
proposed for the proof of his geometrization theorem for 
orbifolds [{\bf T2}]. As a corollary, we will give a proof 
of the following special case of Thurston's geometrization 
theorem.

{\sc Corollary B.} {\em If $M$ is an irreducible, closed, 
atoroidal 3-manifolds, it is not Seifert manifold and
admits a finite group $G$ action. If the order of $G$ is odd,
the $G$-action is effective and not fixed-point-free, then 
the quotient $M/G$ is a geometric orbifold.}

This paper is orgnized in seven sections. \S 1 collects 
some basic facts for geometric cone manifolds and their
limits, we refer reader to somewhere
else for details. Our first goal is using equivariant Ricci
flow to study the 
topology and geometry of compact Euclidean cone manifolds.
A briefly review of Hamilton's work on Ricci flow in given in
\S 2, and we also establish a version of Ricci flow on
orbifold. Using Ricci flow on orbifold we study compact
Euclidean cone manifold in \S 3. Our next goal is a compactness
theorem for hyperbolic cone structures. For the purpose, we
show that the compact hyperbolic cone manifolds with cone
angles less than $2 \pi /3$ cannot become thinner and
thinner everywhere in \S 4 and the compactness theorem will
be given in \S 5. In \S 6, we review the deformation 
theory of hyperbolic cone structure; and the proofs of
Theorem A and Corollary B will be given in the last section.

Essentially this paper is a rewrite of my Ph.D. thesis 
[{\bf Z}]. I would like to express my gratitude 
to my thesis advisor Robert D. Edwards here. Without his
encouragement, this work could not been done. 

\vspace{0.1in}

{\bf \S1 Preliminaries}

We collect some basic facts about geometric cone manifolds 
and their limits in this section,
for details we refer readers to [{\bf Ho}] and [{\bf SOK}].

Let $C$ be a 3-dimensional geometric cone manifold, 
$X_C$ the underlying space, $\Sigma_C$ the singular locus. 
In this paper, we will always assume that $X_C$ is a 
closed manifold, $\Sigma_C$ is a link in $X_C$ and all 
cone angles are less than or equal to $\pi$. The pair 
$(X_C, \Sigma_C)$ is called the combinatorial type of 
the geometric cone manifold $C$. Two cone manifolds $C_1$ 
and $C_2$ are said to be isomorphic if there is an isometry 
between them. The isometry of two geometric cone manifolds
induces an homeomorphism between their combinatorial types.

A geometric cone manifold $C$, each of whose cone angles is of 
the form $2 \pi /n$ for some integer $n \geqslant 2$, is a geometric
orbifold. The concept of the cone manifold is a generalization
of orbifold.

A geometric cone structure is a geometric structure in the
sense of Thurston [{\bf T1}] but in a more general
setting. So we can talk about the developing map
and the holonomy representation. It is clear that 
holonomy representations of two isomorphic cone manifolds $C_1$ and
$C_2$ are conjugate in the sense that the following diagram 
commutes:
$$
\begin{array}{ccc}
\pi_1 (X_1 - \Sigma_1) & \rightarrow & PSL_2(\mathbb{C}) \\
\downarrow & & \downarrow \\
\pi_1 (X_2 - \Sigma_2) & \rightarrow & PSL_2(\mathbb{C})
\end{array}
$$

An $\varepsilon$-ball $B_{\varepsilon}(x)$ in a geometric
cone manifold $C$ is standard if it is standard in the usual sense 
or the ball has a sigular diameter and $x$ is on the sigular 
diameter. For a point $x \in C$, The injectivity radius
of $x$ is defined as follows:

If $x \in \Sigma_C$, 
$$\mbox{inj}(x) = \sup \{\varepsilon |\varepsilon \mbox{-ball }
B_\varepsilon (x)\mbox{ is standard}\};
$$

If $x \not\in \Sigma_C$,
\begin{eqnarray*}
\mbox{inj}(x) & = & \max( \sup \{\varepsilon |\varepsilon \mbox{-ball }
B_\varepsilon (x)\mbox{ is standard}\}, \\
& & \sup \{ \delta | \mbox{there
is } y \in \Sigma_C \mbox{ such that inj}(y) > 2 \delta \mbox{ and } 
d(x, y) < \delta\}).
\end{eqnarray*}

It is easy to see that the injectivity radius is positive lower 
semicontinuous fuction on $C$.

Using injectitivty radius we can divid a geometric cone 
manifold into the thin part and the thick part.
For any $\varepsilon > 0$, the 
$\varepsilon$-thin part of a geometric cone manifold
$C$ is the set $C_{\mbox{\scriptsize{thin,}} \varepsilon}
= \{ x | \mbox{inj}(x)\leqslant \varepsilon \}$ and the 
$\varepsilon$-thick part is the complement of the 
$\varepsilon$-thin part, i.e., 
$C_{\mbox{\scriptsize{thick,}} \varepsilon } = 
C - C_{\mbox{\scriptsize{thin,}} \varepsilon }$.

As the injectivity radius on a Riemannian manifold, 
the injectivity radius on a cone manifold with nonpositive
constant sectional curvature
cannot decrease too fast. More precisely, we have

{\sc Theorem 1.1.}([Proposition 6.1, {\bf SOK}]) 
{\em Let $C$ be a 3-dimensional cone manifold with 
nonpositive constant sectional curvature. 
Suppose that there is a positive constant $\omega$ 
such that all cone angles are at least $\omega$. 
Then, for any $R, \varepsilon > 0$, there is an 
$\delta > 0$ which depends only on 
$R, \varepsilon$, $\omega$ and does not depend on 
$C$, such that, if inj$(x) > \varepsilon$, then 
$B_R(x,C) \subset C_{\mbox{\scriptsize{thick,}} 
\delta}$.}

A useful tool to analyze the geometric cone manifold
is the pointed Gromov limit. For the definition
of Gromov limit and the proofs of the following 
theorems we refer readers to [{\bf G}] and [{\bf SOK}].

{\sc Theorem 1.2.} (see [Proposition 3.6 and Theorem 4.2, 
{\bf SOK}]){\em Let $(C_n,x_n)$ be a sequence of cone manifolds 
with nonpositive constant sectional curvature
and

1) the curvatures $K_n$ has a limit $K_\infty$,

2) all cone angles have a uniform lower bound 
$\omega$, and

3) there is an $\varepsilon >0$, such that 
$x_n\in C_{n\mbox{\scriptsize{,thick,}} \varepsilon }$ for 
all $n$.

\noindent Then it has a convergent subsequence $(C_{n_k}, 
x_{n_k})$, and the limit space $(C,x)$ is a cone manifold 
with constant curvature $K_\infty $.}

In general, even if all the $C_n$ have the same fixed 
combinatorial type, the limit space $C$ does not have to 
have the same combinatorial type. However 
$C$ does have the same topology as $C_n$ locally.

{\sc Theorem 1.3.} (see [Theorem 4.2 and Proposition 8.1,
{\bf SOK}]){\em Let $(C_n,x_n)$ be a convergent sequence of 
cone manifolds with nonpositive constant sectional curvature
such that $\lim (C_n,x_n)=(C,x)$ is still a 
cone manifold. For any $R>0$, if $B_R(x,C)$ is a proper set 
and $n$ large enough, $(B_R(x_n,C_n),B_R(x_n,C_n)\cap \Sigma _n)$ 
is homeomorphic to $(B_R(x,C),B_R(x,C)\cap \Sigma )$ and, 
in fact, this homeomorphism can be chosen as an almost isometry.}

This also shows that the holonomy of $B_R(x_n,C_n)$ can be chosen
to converge to a holonomy of $B_R(x, C)$ if all $C_n$ are hyperbolic
cone manifolds.

To understand the ``local picture'' near a point in a hyperbolic 
cone manifold with a small injectivity radius $\delta$, we can rescale 
the hyperbolic metric by multiplying $\delta^{-1}$, which is a 
large number. Then the new structure has constant curvature
$\delta^2$. If $\delta$ small enough, the new structure
is very close to a Euclidean structure. Theorem 1.3 says that
locally the hyperbolic cone manifold has a topology of 
a Euclidean cone manifold. More precisely, we have the 
following analogue of the Kazhdan-Margulis theorem.

{\sc Theorem 1.4.} ([Proposition 8.1, {\bf SOK}]{\em For 
any $R>0$ and $\omega >0$, 
there is a $\delta < 1 / R$ so that, for any 3-dimensional 
hyperbolic cone manifold $C$ with cone angle at least
$\omega $ and $x\in C_{\mbox{\scriptsize{thin, }}\delta}$, 
we have that $(B_{R\mbox{\scriptsize{inj}}(x)}(x,C),
B_{R\mbox{\scriptsize{inj}}(x)}(x,C)
\cap \Sigma _C)$ is homeomorphic to $(B_R(y,E),B_R(y,E)\cap 
\Sigma _E)$ for some 3-dimensional Euclidean cone manifold 
$E$ and inj$(y)=1$ and after a rescale, the homeomorphism
is an almost isometry.}

This theorem tell us that 3-dimensional Euclidean cone 
manifolds play an important role in analyzing the thin 
part of 3-dimensional hyperbolic cone manifolds. A complete
classification of noncompact Euclidean cone manifolds is
known, and we state the conclusion here and refer readers 
to [\S 5, {\bf SOK}] for the discussion.

{\sc Theorem 1.5.} {\em Let $E$ be a noncompact 3-dimensional 
Euclidean cone manifold, then $E$ is one of the following:

1) $\mathbb{E}^3$ or a product of $\mathbb{E}^1$ and an infinite disk 
with a cone;

2) a product of $\mathbb{E}^1$ with a torus or a compact 2-dimensional
Euclidean cone manifold;

3) a twisted $\mathbb{E}^1$ bundle over a Klein bottle or a projective
plane with two cones of angle $\pi $;

4) a bundle over $S^1$ whose fiber is $\mathbb{E}^2$ or an infinite 
disk with a cone; or

5) a quotient of $\mathbb{E}^3$ modulo one of the following groups:
Denote by $T_{\bf a}$ the translation in a vertor ${\bf a}$ and 
by $R_{l, \varphi}$ the rotation through an angle $\varphi$ about 
an axis $l$,
i) the group generated by $R_{l_1, \pi}, R_{l_2, \pi}$ for two
nonintersect lines $l_1, l_2$; ii) the group generated by $R_{l_1,
\pi}, R_{l_2, \pi}$ and $T_{\bf a}$, where $l_2 = T_{\bf b}l_1$, 
${\bf a}$ is parallel to $l_1$ and perpendicular to ${\bf b}$, 
iii) the group generated by $R_{l_1, \pi}, R_{l_2, \pi}$ and 
$T_{\bf a}$, where $l_1, l_2$ are nonintersect lines, perpendicular 
to each other and ${\bf a}$ is parallel to $l_1$.}

{\sc Remark.} All of the compact 2-dimensional Euclidean cone 
manifolds are easily to be classified. They are a flat torus, 
a double of an acute angled Euclidean triangle, a double of a
rectangle or the boundary 
of a Euclidean tetrahedron with equal opposite edges. For the case
of boundary of a Euclidean tetrahedron, we said that the tetrahedron
must have equal opposite edges is because our convention
that all cone angles are not large than $\pi$, so all four cone
angles are equal $\pi$.  

\vspace{0.1in}

{\bf \S2 The Ricci Flow on Orbifolds}

In this section, we give a brief review of Hamilton's work 
([{\bf Ha1}] and [{\bf Ha2}]) on Ricci flow of 3-dimensional 
nonnegatively curved Riemannian manifolds. As Hamilton pointed 
out that the flow is invariant under an action by isometries, 
we can establish an orbifold version of Ricci flow. 
These results will be used to study the topology and the 
geometry of compact 3-dimensional Euclidean cone manifolds 
in the next section.

Let $M$ be a closed 3-dimensional manifold and $g_0$ a Riemannian 
metric with positive Ricci curvature. Hamilton considered a partial 
differential equation
$$
\left\{ 
\begin{array}{l}
\partial_t g = \frac{2}{3} rg - 2 \mbox{Ric} \\
g(\cdot ,0)=g_0 
\end{array}
\right. 
$$
on $M$, where $r=\int R/\int 1$ is the average of the scalar 
curvature $R$. Hamilton, using the Nash-Morse inverse function 
theorem, showed that the equation has a short time solution for 
any initial metric $g_0$ and then estimated the curvatures by 
the maximum principle for parabolic equations under the assumption 
that the initial metric $g_0$ has positive Ricci curvature. 
The main result in [{\bf Ha1}] is the following theorem.

{\sc Theorem 2.1.} {\em For a compact 3-dimensional manifold $M$ 
with positive Ricci curvature, the metric evolution
$$
\partial_t g= \frac{2}{3} rg - 2 \mbox{Ric} 
$$
has a solution for all the time and converges to a metric of 
positive constant Riemannian curvature as $t \rightarrow \infty $.}

The equation is invariant under the full diffeomorphism group of $M$, 
so any isometry for the initial metric are preserved as the metric 
evolves. This fact allows us to establish an orbifold version of 
Hamilton's theorem.

{\sc Definition 2.2.} Let $O$ be an orbifold and $g$ a Riemannian
metric on the complement of the singular locus, we say that it is a
Riemannian metric on the orbifold $O$ if, passing to a local manifold
cover, the lift of $g$ can be extended to a smooth Riemannian metric. 
For a Riemannian orbifold $O$, we can talk about the connection and the
curvature.

{\sc Theorem 2.3.} {\em Suppose that a compact very good orbifold $O$ 
admits a Riemannian metric with positive Ricci curvature, then it is 
a spherical orbifold.}

An orbifold is said to be very good, it means that the orbifold 
is the quotient of a manifold modulo a finite group. 

{\sc Proof.} Since $O$ is very good, $O$ is a quotient of a 
manifold $M$ modulo a finite group $G$. We lift the Riemannian metric on 
$O$ to get a metric $g_0$ with positive Ricci curvature on $M$.
Obviously, $G$ consists of isometries of $g_0$. Applying Hamilton's 
theorem to $(M,g_0)$, we have a Ricci flow $g_t$, which converges 
to a metric $g_\infty $ of a positive constant Riemannian curvature 
as $t\rightarrow \infty $, and $G$ is a group of isometries on $g_t$ 
for all $t$ and $g_\infty$. So $(M,g_\infty )/G$ is a spherical 
structure on the orbifold $O$. \hspace*{\fill}$\square $

Hamilton also developed a method to deal with the case of nonnegatively
curved manifold in [{\bf Ha2}]. He considered the unnormalized equation
$$
\partial_t g= -2 \mbox{Ric}. 
$$
The equation has a short time solution for any initial metric $g_0$. 
This equation for the metric implies a heat equation for the Riemannian 
curvature tensor $R_{ijkl}$:
\begin{eqnarray*}
\partial_t R_{ijkl} & = & \Delta
R_{ijkl}+2(B_{ijkl}-B_{ijlk}+B_{ikjl}-B_{iljk})\\ 
& & -g^{pq}(R_{pjkl}R_{qi}+R_{ipkl}R_{qj}+R_{ijpl}R_{qk}
+R_{ijkp}R_{ql}),
\end{eqnarray*}
where $B_{ijkl}=g^{pr}g^{qs}R_{ipjq}R_{krls}$.

After choosing a family of isometries $u_t:(V,h)\rightarrow (TM,g_t)$, 
Hamilton pullbacked the Levi-Civita connections and the 
curvatures on $(TM,g_t)$ and defined the covariant derivatives and 
the Laplacian on $(V,h)$. Then the pullback of the heat equation for 
Riemannian curvature tensor can be written simply as
$$
\partial_t R_{abcd}=\Delta
R_{abcd}+2(B_{abcd}-B_{abdc}+B_{acbd}-B_{adbc}). 
$$

The Riemannian curvature tensor can be regarded as a symmetric 
bilinear form $K$, which is called the curvature operator, on 
the 2-forms $\Lambda ^2(V)$. For a 2-form $\alpha \wedge \beta$,
$K( \alpha \wedge \beta)$ is the product of $| \alpha \wedge \beta|$
and the sectional curvature of the plane with the Pl\"ucker 
coordinates $\alpha\wedge \beta$. The curvature operator contains the
information as much as the Riemannian curvature tensor and the 
positive curvature operator implies positive sectional curvature.
However in general positive curature operator is not equivalent
to positive sectional curvature. 
In dimension 3, any 2-form is pure, which means that
any 2-form can be written as a wedge product of two 1-forms. 
So the positive curvature operator and the positive sectional 
curvature are the same thing in this dimension. Now the heat 
equation for the Riemannian curvature tensor can be rewritten 
as an equation for the curvature operator
$$
\partial_t K = \Delta K + K^2 + 2 \mbox{Adj}(K). 
$$
The crucial lemma in [{\bf Ha2}] is the following.

{\sc Lemma 2.4.} [Lemma 8.2., {\bf Ha2}] {\em Let $Q$ be a symmetric
bilinear form on $V$. Suppose $Q$ satisfies a heat equation $\partial_t
Q = \Delta Q + \phi (Q)$, where the matrix $\phi (Q)\geqslant 0$ for all 
$Q \geqslant 0$. Then if $Q \geqslant 0$ at time $t=0$, it remains so for 
$t\geqslant 0$. Moreover there exists an interval $0<t<\delta $ on 
which the rank of $Q$ is a constant, the null space of $Q$ 
is invariant under parallel translation and invariant in the 
time and also lies in the null space of $\phi (Q)$.}

Suppose that $M$ is a 3-dimensional Riemannian manifold with 
nonnegative curvature operator, we solved the unnormalized Ricci flow
$$
\partial_t g = - 2 \mbox{Ric}. 
$$
Applying Lemma 2.4 to the evolving equation for the curvature 
operator, we can see that, if we start with a nonnegative 
curvature operator, after the flow proceeds for a while we 
will get $K=0$, $K>0$ or that the rank of $K$
is constant one or two. If $K=0$, then we have a Euclidean 
structure which is the same as the metric we started with. 
If $K>0$, we can use Theorem 2.1 to deform the metric to a 
spherical structure further. If the rank of $K$ is one, then 
the image of $K$, which is the orthogonal complement of the null
space, has dimension $1$ and also is invariant under the parallel 
translation. This one dimensional space is generated by a 
2-form $\varphi $. We know that $\varphi $ is a pure element, 
i.e., $\varphi $ is the Pl\"ucker coordinates of a plane field 
in $TM$, which is invariant under parallel translation and
the sectional curvature of each plane is the positive eigenvalue of the
curvature operator. This gives an orthogonal decomposition 
$TM=V_1\oplus V_2$, with $V_i$ being invariant under parallel 
translation. The universal covering space $\tilde M$ also has 
a such decomposition. Applying the de Rham decomposition 
theorem to the universal cover $\tilde M$, $\tilde M$ splits isometrically 
into a product $N\times \mathbb{E}^1$, where $N$ is a
positively curved surface and the curvature has a positive 
lower bound since $M$ is compact. This implies that $N$ is a sphere. 
The metric on $N$ may not be standard, but we will show that 
we can replace this metric by a standard one such that any given
group of isometries of the original metric on 
$N\times \mathbb{E}^1$ remains a group of isometries of 
$S^2\times \mathbb{E}^1$.

It is clear that isometry group acting on $N\times \mathbb{E}^1$ 
splits into a product of two isometric actions of $N$ and 
$\mathbb{E}^1$. The action on $N$ induces a map into Conf$(S^2)$. 
If this action can be conjugated to a subgroup of $O(3)$, 
then we can replace the metric on $N$ by a standard
one such that the given group of isometries of the original metric on 
$N\times \mathbb{E}^1$ remains a group of isometries of $S^2\times 
\mathbb{E}^1$. To do this, we need a lemma to
characterize those subgroup of Conf$(S^2)$ which can 
be conjugated into $O(3) $. Note that the Conf$(S^2)$ 
action on $S^2$ has a natural extension to
an action on $\mathbb{H}^3$, and a subgroup of Conf$(S^2)$ can 
be conjugated into $O(3)$ if and only if the action has 
a common fixed point in $\mathbb{H}^3$.

{\sc Lemma 2.5.} {\em A subgroup $\Gamma \subset $Conf$(S^2)$ 
fixes a point in $\mathbb{H}^3$ if and only if $\Gamma \cap 
PSL_2(\mathbb{C})$ consists of only elliptic elements (and $I$).}

{\sc Proof.} It is clear that if $\Gamma $ fixes a point in 
$\mathbb{H}^3$ then $\Gamma \cap PSL_2(\mathbb{C})$ consists of 
only elliptic elements. On the other hand, it is known that 
if $\Gamma \cap PSL_2(\mathbb{C})$ consists of only elliptic 
elements, $\Gamma \cap $$PSL_2(\mathbb{C})$ fixes a point in $\mathbb{H}^3$ 
(see [Theorem 4.3.7, 
{\bf B}]). So if $\Gamma \subset PSL_2(\mathbb{C})$, then we are done.
Otherwise, $\Gamma _1=\Gamma \cap PSL_2(\mathbb{C})$ is a normal 
subgroup in $\Gamma $ with index $2$. Fix a $g\in \Gamma -
\Gamma _1$, we have $g($Fix$(\Gamma _1))=$Fix$(\Gamma _1)$ where 
Fix$(\Gamma _1)$ is the fixed point set of $\Gamma _1$ in 
$\mathbb{H}^3$, which is $\mathbb{H}^3$, a geodesic or a point.
If Fix$(\Gamma _1)=\mathbb{H}^3$, then $\Gamma $ is a group of order $2$,
consists of $I$ and a reflection $g$. Of course $\Gamma $ has 
a fixed point in $\mathbb{H}^3$. If Fix$(\Gamma _1)$ is a geodesic 
and $g$ have no fixed point on it, then $g^2\in \Gamma _1$ also 
has no fixed point on it. This is impossible, so $g$ must fixes 
a point on Fix$(\Gamma _1)$. If Fix$(\Gamma_1) $ is only a point, 
then $g$ fixes this point. This shows that $\Gamma
=\Gamma _1\cup g\Gamma _1$ fixes this point in $\mathbb{H}^3$. 
\hspace*{\fill} $\square $

Now we only need to show that any orientation preserving element 
in a group of isometries of $N$ must be elliptic. If $g$ 
is parabolic, then $g$ has a unique fixed point $p\in N$. For 
any point $q\in N - \{p\}$, $g^n(q)\rightarrow p$, this 
contradicts the hypothesis that $g$ is an isometry on $N$. 
A similar argument shows that there is no
hyperbolic element in a group of orientation preserving 
isometries of $N$. Now Lemma 2.5 says that we can replace 
the original metric on $N$ by a standard one such that any 
given group of isometries of $N\times \mathbb{E}^1$ remains a
group of isometries of $S^2\times \mathbb{E}^1$.

Hamilton shown that the case that rank of $K$ is two cannot 
really happen. Diagonalize $K$, then $K^2$ and Adj$(K)$ are 
also diagonalized, say
$$
K=\left( 
\begin{array}{ccc}
\lambda &  &  \\  
& \mu &  \\  
&  & \nu 
\end{array}
\right) ,\mbox{ then }K^2=\left( 
\begin{array}{ccc}
\lambda ^2 &  &  \\  
& \mu ^2 &  \\  
&  & \nu ^2 
\end{array}
\right) ,\mbox{ and Adj}(K)=\left( 
\begin{array}{ccc}
\mu \nu &  &  \\  
& \nu \lambda &  \\  
&  & \lambda \mu 
\end{array}
\right) . 
$$
If the rank of $K$ is two, then we can assume that $\lambda =0$ 
but $\mu \nu \neq 0$, Lemma 2.4 says that the null space of $K$ 
lies in the null space of $K^2+2$Adj$(K)$, in other words, 
$\mu \nu =0$. We get a contradiction.

The equation
$$
\partial_t g = -2 \mbox{Ric} 
$$
is also invariant under the full diffeomorphism group of $M$, so we can
also establish an orbifold version. To summarize, we have

{\sc Theorem 2.6.} {\em Let $O$ be a very good orbifold which admits a
Riemannian metric with nonnegative curvature operator. Then it is a
geometric orbifold locally modelled in $S^3$, $S^2\times \mathbb{E}^1$ 
or $\mathbb{E}^3$. Furthermore, if $O$ has Euclidean structure, then 
the original Riemannian metric must be Euclidean.}

{\sc Remark.} Hamilton actually did prove a stronger result that, 
if $M$ has nonnegative Ricci curvature, then $M$ admits a 
Euclidean geometry, a spherical geometry or a geometry locally 
modelled in $S^2\times \mathbb{E}^1$. We can get an orbifold 
version for this stronger result in the same
manner, but Theorem 2.6 is enough for our purpose.

\vspace{0.1in}

{\bf \S3 Compact Euclidean Cone Manifolds}

Using the Ricci flow on orbifolds, we will discuss the topological type of 
compact 3-dimensional Euclidean cone manifolds in this section, 
and show that at least one cone angle in a compact Euclidean cone 
manifolds is not less than $2\pi /3$.

Suppose $E$ is a compact 3-dimensional Euclidean cone manifold,
$X_E$ be the underlying space, and $\Sigma_E$ be the singular 
locus. If $\Sigma_E $ is empty or all cone angles are $\pi $, 
$E$ is a Euclidean manifold or orbifold, they are a quotient 
of $\mathbb{E}^3 $ by a 3-dimensional crystallographic group. 
All 3-dimensional crystallographic groups have been classified 
(see for example [{\bf J}]). So, without loss of generality, 
we can assume that $\Sigma _E\neq \varnothing $ and that at least 
one cone angle is strictly less than $\pi $. 

Near the singular locus, we use Fermi coordinates
$$
ds^2 = dt^2 + dr^2 + r^2 \alpha^2 d\theta^2 /4 \pi^2,
$$
where $\alpha$ is the cone angle. Geometrically, it is a product 
of $\mathbb{E}^1$ and a disk with a cone. 
For another metric in a cylidrical coordinates system in the form
$$
ds^2 = dt^2 + dr_1^2 + f^2(r_1)d\theta^2,
$$
the curvature operator will be
$$
\left( \begin{array}{ccc}
0 & & \\
& 0 & \\
& & -f''(r_1)/f(r_1)
\end{array} \right).
$$
If near $0$, $f(r_1) = \beta \sin r_1$ and $f(r_1) > 0$ 
for all $r_1 >0$, this is a singular
Riemannian metric with cone angle $2 \pi \beta$ at $r_1 = 0$.
The metric is smooth at $r_1=0$ if $\beta =1$. Pick a function 
$f(r_1)$ such that $f(r_1) = \sin r_1$ for $r_1$ near $0$, 
$f''(r_1) \leqslant 0$ and 
$$
f(r_1) = (r_1 + (2 \pi - \alpha) \varepsilon / 2 \alpha ) 
\alpha / 2\pi \mbox{ for } r_1 > \varepsilon > 0.
$$
This can be done because $\alpha < 2 \pi$ and in this case 
$f''(r_1)$ cannot be identically zero. This is a metric with 
nonnegative curvature operator and, in the region $r_1 > 
\varepsilon$, the metric is globally isometric to $ds^2 = dt^2 
+ dr^2 + r^2 \alpha^2 d \theta^2 / 4 \pi^2$ under the coordinates 
change $r = r_1 + (2 \pi - \alpha)\varepsilon / 2 \alpha$.
So we get a metric with nonnegative curvature operator on the 
underlying space $X_E$. Note that this metric is not flat. Now 
we can apply Theorem 2.6 to conclude that $X_E$ must be a 
spherical manifold or covered by $S^2\times \mathbb{E}^1$.

{\sc Remark.}  Acturally we can increase the cone angle to any 
amount not necessary to $2\pi $ by deforming the Euclidean cone
structure to a metric with nonnegative curature operator in the same
fashion. If all cone angles are in the form of $2\pi /n$, we get 
a Riemannian orbifold with nonnegative curvature operator.

If $X_E$ is covered by $S^2\times \mathbb{E}^1$, then we can pullback 
the singular Euclidean structure on $X_E$ to a one on the universal 
cover $S^2\times \mathbb{E}^1$. Denote this Euclidean cone structure 
by $\tilde E$. Since $\tilde E$ has two ends and at least one cone angle
is less than $\pi$, by the classification, 
$\tilde E$ is a product of $\mathbb{E}^1$ and a double of an acute 
angled Euclidean triangle. There are two possible deck transformation 
groups $\mathbb{Z}$ or $\mathbb{Z}_2*\mathbb{Z}_2$. Suppose that the deck group 
is $\mathbb{Z}$, the $\mathbb{Z}$ acts on $\mathbb{E}^1$ as a translation 
group. This shows that $E$ can be regarded as a quotient of
a product of an interval $I$ and a double of a Euclidean triangle 
with an isometric mondromy identifying the top and the bottom. If 
the deck group is $\mathbb{Z}_2*\mathbb{Z}_2$, then it must act on 
$\mathbb{E}^1$ as an isometry group generated by two reflections. 
Since some elements have a fixed point on $E^1$, these elements 
will keep a cross intersection invariant. On the other hand, a 
double of a Euclidean triangle does not allow an orientation 
reversing isometric inversion which permute three cone points, 
so the cross intersection cannot be invariant,
i.e., this is impossible. In this case, $X_E$ must be $S^2\times S^1$, 
$\Sigma _E$ can be any 3-string braid in $S^2\times S^1$ and at least 
one cone angle is not less than $2\pi /3$. Note that, in this case, 
$X_E-\Sigma_E$ is Seifert fibered and the following theorem 
shows that the converse is also true.

{\sc Theorem 3.1.} {\em Let $(X,\Sigma )$ be a combinatorial type of a
compact Euclidean cone manifold $E$ with nonempty singular locus and at
least one cone angle less than $\pi $. Then $N=X-{\mathcal N}(\Sigma )$ 
is an irreducible, atoroidal 3-manifold with an incompressible 
boundary, where ${\mathcal N}(\Sigma )$ is a regular neighborhood of 
$\Sigma $.}

{\em Furthermore, $N$ is Seifert fibered if and only if $X$ is 
$S^1\times S^2 $.}

{\sc Proof.} It is not difficult to see that any essential simple 
closed loop in $\partial N$ will be sent to a nontrivial element 
under the holonomy, so $\partial N$ is incompressible.

To show that $N$ is irreducible and atoroidal, we will use minimal 
surface techniques. We deform the singular Euclidean structure near 
the singular locus $\Sigma $ as before. First we can find a number 
$\varepsilon >0$ such that $B_\varepsilon (\Sigma )$ is a disjoint 
union of $B_\varepsilon (C_i)$ for all component $C_i$ of $\Sigma $ 
and each $B_\varepsilon (C_i)$ is a twisted product of $S^1$ and a 
disk with a cone. Without loss of generality we can assume 
that $N=X-B_{\varepsilon /10}(\Sigma )$. 

Replace the metric on $B_{\varepsilon}(\Sigma) - B_{\varepsilon 
/10}(\Sigma)$ by
$$
ds^2 = dt^2 + dr^2 + f^2(r) d\theta^2,
$$
where $f''(r) \geqslant 0$, $f(r) = \alpha r /2 \pi$ when $r > \varepsilon$
and $f(r) =$ constant when $r$ near $\varepsilon /10$.
Thus we have a Riemannian metric with nonpositive 
sectional curvature on $N$ such that the boundary is totally
geodesic.

If $N$ is reducible, we can find a least area immersed sphere $S$. 
The sectional curvature is nonpositive and, since $S$ is minimal, 
the determinant of the second fundamental form on $S$ is also 
nonpositive. So the Gaussian curvature of this immersed sphere is 
nonpositive, because it is the sum of the sectional curvature of 
the tangent plane and the determinant of the second fundamental form. 
On the other hand, $\int_SKd\sigma =4\pi >0$ by Gauss-Bonnet theorem. 
This is a contradiction and it shows that $N$ is irreducible.

If we start with an incompressible torus in $N$, we can also find an
immersed least area torus $T$ in the same homotopy class. In fact this $T$
is either embedded or a 2-fold cover onto an embedded Klein bottle $K$. By
the same argument as above, we know that the Gaussian curvature of this
torus is nonpositive and the Gauss-Bonnet theorem implies that the Gaussian
curvature of this torus will be identically equal to $0$ and, therefore this 
$T$ is totally geodesic and all sectional curvatures on the tangent planes
of $T$ is zero. Now it is easy to see that $T$ is either $\partial N$ 
or contained in $X-B_{\varepsilon}(\Sigma )$. We will argue that $T$ 
is $\partial N$. 
If $T$ is contained in $X-B_{\varepsilon}(\Sigma )$, we get a 
totally geodesic embedded torus $T$ or a Klein bottle $K$ in $E$. 
Cutting $E$ along $T$ or $K$, we can get an $E^{\prime }$ which 
is a compact Euclidean cone manifold with one or two totally geodesic 
flat tori as boundary. Then we can extend $E^{\prime }$ to a 
noncompact Euclidean cone manifold with one or two ends look 
like a product of $\mathbb{E}_{+}^1$ and a torus. On the other hand,
we assume that there is at least one cone angle is less than $\pi $. 
This is a contradiction, since there is no noncompact 
Euclidean cone manifolds with ends look like a product of $\mathbb{E}_{+}^1$
and a torus, and a cone angle less than $\pi $.

To prove the second statement, we look at the holonomy $\rho :\pi
_1(N)\rightarrow $Isom$^{+}(\mathbb{E}^3)$ of $E$. There are three
possibilities for a nontrivial element $\alpha $ in Isom$^{+}(\mathbb{E}^3)$:

a) $\alpha $ has a unique invariant line;

b) $\alpha $ is a $\pi$ rotation along an axis and this axis is the
unique line which is fixed by $\alpha $; or

c) $\alpha $ is a pure translation in a direction $n$, in this case, 
$\alpha $ has infinitely many invariant lines and they are all 
parallel to the direction $n$.

Suppose that $N$ is Seifert fibered. Then $\pi _1(N)$ has a 
central infinite cyclic subgroup $\mathbb{Z}$ which is generated by 
an element $a$ represented by a regular fiber. We can choose this 
regular fiber on $\partial N$, so $\rho (a)$ is nontrivial in 
Isom$^{+}(\mathbb{E}^3)$. If $\rho (a)$ is type a) or b), there is 
a unique line $l_a$ which is invariant or fixed under $\rho(a)$. 
$a$ is a central element, so for any element $b$, $l_a$ is also
invariant or fixed under $\rho (b)$. This is impossible, since 
we cannot get a compact Euclidean cone manifold in such way.

Now we know that $\rho (a)$ must be a pure translation along the direction 
$n $. For any other element $b$, $\rho (b)$ commutes with $\rho (a)$. 
If $\rho (b)$ is type a) or b), the unique invariant line or fixed line 
will be parallel to the direction $n$. This shows that the foliation 
which arises by all planes perpendicular to $n$ is invariant under 
$\rho (\pi _1(N))$, so it gives us a foliation on $E$. Each leaf 
of this foliation must be a compact 2-dimensional Euclidean cone 
manifold, otherwise we will have a noncompact 2-dimensional cone 
manifold with infinitely many cone points. This is impossible, since 
the only noncompact 2-dimensional cone manifold with cone
angles less than $\pi $ is an infinite disk with a cone. This foliation
gives us an $S^2$ bundle structure on $X$, so $X$ must be $S^1\times S^2$. 
\hspace*{\fill}$\square $

Finally, we show that there is no compact Euclidean cone manifold with all
cone angles less than $2\pi /3$.

{\sc Theorem 3.2.} {\em For any compact Euclidean cone manifold, there is at
least one cone angle not less than $2\pi /3$.}

{\sc Proof.} We have already known that at least one cone angle must be not
less than $2\pi /3$, if the underlying space is $S^1\times S^2$. If the 
theroem does not hold, then we can assume that there is a Euclidean 
cone manifold $E$ with all cone angles less than $2\pi /3$ and $X_E=S^3$. 
This is because any Euclidean cone structure on a spherical manifold can 
be lifted to a structure on $S^3$.

Now we increase all cone angles to $2\pi /3$ to get a nonnegatively curved
Riemannian orbifold $O$ as we did at the beginning of the section. $O$ can 
be triply covered by a manifold $N$, since such an $N$ can be
constructed by using a Seifert surface of $\Sigma _E$. By Theorem 1.6, now
we can conclude that $O$ is a geometric orbifold with spherical
geometry or geometry of $S^2\times \mathbb{E}^1$. In fact, $O$ cannot
admit $S^2\times \mathbb{E}^1$ geometry. Otherwise, $N$ will be $RP^3\sharp
RP^3 $ or $S^2\times S^1$. By the equivariant minimal sphere theorem, the
quotient of $RP^3\sharp RP^3$ or $S^2\times S^1$ modulo a $\mathbb{Z}_3$ action
will never yield an $S^3$. So, in this case, $O$ is a spherical
orbifold. The Riemannian universal cover of $N$ is $S^3$ and, if we 
lift the identity on $O$, we get the fundamental group 
$\pi^{\mbox{\scriptsize{orb}}}(O)=G$ of $O$ which is a subgroup of $SO(4)$ 
and $O =S^3/G$.

To analyze this group $G$, we regard $SO(4)$ as $S^3\times S^3/\mathbb{Z}_2$.
There is a natural map $p:SO(4)\rightarrow SO(3)\times SO(3)$. Let $H=p(G)$
and let $H_1$, $H_2$ be the projections of $H$ into the two factors of 
$SO(3)\times SO(3)$. Note that all nontrivial elements in $G$ which have a
fixed point in $S^3$ have order three, this allows us to establish the
following lemma.

{\sc Lemma 3.3.} {\em $H_1$ or $H_2$ is cyclic.}

{\sc Proof of the lemma.} If the order of $G$ is odd, both of $H_1$ and 
$H_2$ have odd order. All finite subgroups of $SO(3)$ with odd order are 
cyclic, so the lemma follows.

We suppose that $G$ has even order and then $G$ must have an element of
order two. By our assumption that all cone angles in $O$ is $2 \pi /3$,
any elements in $G$ of order two are fixed point free and the
only fixed point free involution in $SO(4)$ is $-I$ which lies in the 
kernel of $p$, so we have $G=p^{-1}(H)$. Therefore, any nontrivial 
elements which act on $SO(3)$ with a fixed point must have order three 
as well, because we can lift a such element to one in $G$ which has a 
fixed point in $S^3$. In fact, we only need that any element of order 
two acts on $SO(3)$ freely. Another fact we need is that $(u_1,u_2)
\in SO(3)\times SO(3)$ acts on $SO(3)$ freely if and only if $u_1$ 
is not conjugate to $u_2$ in $SO(3)$. In particular, $H$ cannot 
contain an element $(u_1,u_2)$ where both $u_1$ and $u_2$ have order 
two, since all elements of order two are fixed point free and therefore 
they are conjugate in $SO(3)$.

We can assume that $H$ has even order, otherwise, both $H_1$ and 
$H_2$ have odd order and, therefore, they are cyclic. Let 
$H_i^{\prime }=H\cap H_i$, then $H_i^{\prime }$ is normal in 
$H_i$ and $H_1/H_1^{\prime }=H/(H_1\times H_2)=H_2/H_2^{\prime }$.

Since $H_1^{\prime }\times H_2^{\prime }$ does not contain an element 
$(u_1,u_2)$ with both $u_1$ and $u_2$ have order two, one of 
$H_1^{\prime }$ or $H_2^{\prime }$ must have odd order. Say 
$H_1^{\prime }$ has odd order, so it must be a cyclic group. 
We also want to show that $H_2^{\prime }$ contains all elements of 
order two in $H_2$. Suppose that $u_2\in H_2$ has
order two. We can find an element $u_1\in H_1$ such that $(u_1,u_2)\in H$. 
$(u_1,u_2)^2=(u_1^2,I)$ implies that $u_1^2\in H_1^{\prime }$, so $u_1^2$ has
odd order $n$. Hence either $u_1^n$ is trivial or has order two. On the
other hand, $(u_1,u_2)^n=(u_1^n,u_2)$, $u_1^n$ cannot have order two. So 
$u_1^n=I$ and $u_2\in H_2^{\prime }$.

Suppose that $H_2$ is not cyclic, $H_2$ must be a dihedral group, $A_4$, 
$S_4 $ or $A_5$. If $H_2$ is not an $A_4$, $H_2=H_2^{\prime }$, since those
groups are generated by elements of order two. Then $H_1=H_1^{\prime }$ is
also cyclic. If $H_2$ is an $A_4$, either $H_2=H_2^{\prime }$ which implies 
$H_1=H_1^{\prime }$, or $H_2^{\prime }=\mathbb{Z}_2\times \mathbb{Z}_2$ 
and, then $H_2/H_2^{\prime }$ has order three and $H_1$ also has odd 
order. So we finish the proof.\hspace*{\fill}$\square $

Now we return to the proof of Theorem 3.2. We can suppose that $H_1$ is
cyclic and then $G$ can be conjugated in $SO(4)$ to a subgroup of $S^1\times
S^3/\mathbb{Z}_2$. The group $S^1\times S^3/\mathbb{Z}_2$ preserves the Hopf
fiberation on $S^3$ and the orientation on fibers, so does $G$. This shows
that $S^3/G$ is Seifert fibered and the singular locus is a union of fibers.
It is in contradiction with Theorem 3.1, which says that $X_E-\Sigma _E$ is
not Seifert fibered. So we cannot have a Euclidean cone manifold with all
cone angles less than $2\pi /3$. \hspace*{\fill}$\square $

Finally we want to point out that
Theorem 3.2 is sharp. $(S^3, \mbox{figure eight})$ does have a 
Euclidean cone structure with cone angles equal to $2\pi /3$.

\vspace{0.1in}

{\bf \S 4 The Foliation on Thin Parts}

In this section, we will show that the injectivity radius cannot be 
very small everywhere on a compact 3-dimensional hyperbolic cone 
manifold.

{\sc Theorem 4.1.} {\em Suppose that $\Sigma $ is a link in $X$ 
and all ends of $X$ are cusps. If $(X,\Sigma )$ is a combinatorial 
type of a 3-dimensional hyperbolic cone manifold with finite volume, 
then $X-\Sigma $ supports a complete hyperbolic structure with 
finite volume.}

{\sc Remark.} In fact, the condition that $X-\Sigma $ is hyperbolic 
is also sufficient for the pair $(X,\Sigma )$ to be a combinatorial 
type of a 3-dimensional hyperbolic cone manifold. We will see that
in \S 6.

{\sc Proof.} To prove the theorem, it is enough to show that the compact
core $N$ of $X-\Sigma $ is an irreducible, atoroidal 3-manifold with
incompressible boundary and it is not Seifert fibered.

Since any simple essential closed loop on the $\partial N$ will be 
sent to a nontrivial element in $PSL_2(\mathbb{C})$ under holonomy, $N$ 
must have an incompressible boundary. Now suppose that $N$ is 
Seifert fibered. Then $\pi_1(N)$ has a central infinite cyclic group 
$\mathbb{Z}$ and the holonomy image of the generator of this central 
cyclic group is nontrivial, since the
generator is represented by an essential loop on the $\partial N$. The fact
that any element in $\pi _1(N)$ commutes with this generator implies that
the holonomy image of $\pi _1(N)$ has a common fixed point in $S^\infty $.
This is impossible, since we cannot get a cone manifold $C$ with finite
volume in a such way.

To prove that $N$ is irreducible and atoroidal, we use an argument similar
to the one used in the proof of Theorem 3.1. First we deform the hyperbolic
structure on $C$ to a Riemannian metric with nonpositive sectional curvature
on $N$, such that, away from a small regular neighborhood of the boundary,
we have constant sectional curvature $-1$; and near the boundary, the metric
is flat and the boundary is totally geodesic.

Near a singular component, we use Fermi coordinates
$$
ds^2=\cosh {}^2rdt^2+dr^2+\sinh {}^2r\alpha ^2d\theta ^2/4\pi ^2. 
$$
In general, if a metric in a cylindrical coordinates system is given by 
$$
ds^2=g^2(r)dt^2+dr^2+f^2(r)d\theta ^2, 
$$
then the curvature operator will be
$$
\left( 
\begin{array}{ccc}
-g^{\prime \prime }(r)/g(r) &  &  \\  
& -g^{\prime }(r)f^{\prime }(r)/f(r)g(r) &  \\  
&  & -f^{\prime \prime }(r)/f(r) 
\end{array}
\right) . 
$$
Now we can choose positive functions $f(r)$ and $g(r)$ nondecreasing and 
concave upward, such that, near $r=0$, both of them are constant and, 
away from $0$, $f(r)$ coincides with $\sinh {}^2r\alpha ^2/4\pi ^2$, 
$g(r)$ coincides with $\cosh {}^2r$.

For a cusp end, we have the metric
$$
ds^2=\frac{dx^2+dy^2+dz^2}{z^2},\quad z>z_0. 
$$
We will deform this metric to a new metric which has the following form:
$$
ds^2=\frac{dx^2+dy^2+dz^2}{f^2(z)}. 
$$
The curvature operator of this metric is given by
$$
\left( 
\begin{array}{ccc}
-f^{\prime 2}(z) &  &  \\  
& f^{\prime \prime }(z)f(z)-f^{\prime 2}(z) &  \\  
&  & f^{\prime \prime }(z)f(z)-f^{\prime 2}(z) 
\end{array}
\right) . 
$$
Now the question of finding the desired metric becomes to find a positive
function $f(z)$ such that $f(z)=z$ near $z=z_0$, $f(z)$ is constant for $
z\gg z_0$, and $f^{\prime \prime }(z)f(z)-f^{\prime 2}(z)\leqslant 0$. Note that $
f^{\prime \prime }(z)f(z)-f^{\prime 2}(z)$ is the numerator in $(f^{\prime
}(z)/f(z))^{\prime } = (f''(z) f(z) - f^{\prime 2}(z))/f^2(z)$. Choose a function $g(z)$ so that $g(z)=1/z$ near $
z=z_0$, $g(z)=0$ for $z\gg z_0$ and $g(z)$ is nonincreasing. Solve $
(f^{\prime }(z)/f(z))^{\prime }=g(z)$ we can get the desired $f(z)$.

Thus we get a metric on $N$ as expected. The irreducibility of $N$ is
followed from the fact that there is no minimal immersed sphere in $N$. If
we start with an incompressible torus in $N$, we will end up an immersed
least area torus $T$ in the same homotopy class. This $T$ is either
embedded or a 2-fold cover to an embedded Klein bottle $K$. By the same
argument as before, this can happen only if it torus lies in a regular 
neighborhood of $\partial N$, and then we can homotopy this torus into 
$\partial N$. We complete the proof of the theorem.\hspace*{\fill}$\square $

{\sc Theorem 4.2.} {\em For any $\omega <2\pi /3$, 
there is a universal constant $\delta >0$, so that, for any 
3-dimensional hyperbolic cone manifold $C$ with all cone angles 
less than $\omega $, we have $C_{
\mbox{\scriptsize{thick, }}\delta }\neq \varnothing $.}

{\sc Proof.} Suppose the contrary that we have a compact hyperbolic 
cone manifold $C$ with all cone angles between $\alpha $ and 
$\omega $, so that $C_{\mbox{\scriptsize{thick},}\delta }=\varnothing $, 
where the $\delta$ is the constant in Theorem 1.4 for $R=120$ 
(later we will see why pick $R=120$ here). 
By Theorem 1.4, for any point $x\in C$, $(B_{120\mbox{\scriptsize{inj}}
(x)}(x,C),B_{120\mbox{\scriptsize{inj}}(x)}(x,C)\cap \Sigma _C)$ 
is homeomorphic to $(B_{120}(y,E),B_{120}(y,E)\cap \Sigma _E)$ for a 
Euclidean 
cone manifold $E$, $y \in E$ and inj$(y)=1$. Recall the proof of 
Theorem 4.1, these Euclidean cone manifolds are limits of the 
rescaled hyperbolic cone manifolds. The hyperbolic cone
manifolds considered here all have cone angles less than 
$\omega <2\pi /3$, and Theorem 3.2 says that there is no 
compact Euclidean cone manifolds with all cone angles less than 
$2\pi /3$, so we can assume further that these
Euclidean cone manifolds which can occur in our case are all noncompact. 
Now, by the classification of noncompact Euclidean cone manifolds, 
we can conclude that they are a product of $\mathbb{E}^1$ 
with a torus, a twisted $\mathbb{E}^1$ bundle over a Klein bottle or 
a bundle over $S^1$ whose fiber is $\mathbb{E}^2$ or an infinite disk 
with a cone. The case that $E$ is a twisted $\mathbb{E}^1$ bundle over 
a Klein bottle will not happen. Otherwise, we could embed this Klein 
bottle in $X-\Sigma$. $X-\Sigma $ is atoroidal and the boundary of
the regular neighborhood of this Klein bottle is incompressible inward, 
so this torus must be parallel to a singular circle in $\Sigma$. 
This contradicts the fact that $X-\Sigma$ is not Seifert fibered.

Now we take a closed look at Euclidean case first, the picture in
this case is the same when we look at a point in $C$ which has
a small injectivity radius.

Case 1). $E$ is a product of $\mathbb{E}^1$ with a trous. The isometry 
group acts on $E$ transitively, so we only need to work at one point 
$y \in E$.
Injectivity radius of $y$ in $E$ is the same as the injectivity radius
of $y$ in the cross section, $B_{\mbox{\scriptsize{inj}}(y)}(y, E)$
is ball with one, two or three pairs of points identified on the
boundary. The image of $\pi_1 (B_{\mbox{\scriptsize{inj}}(y)}(y,E))$
in $\pi_1 (B_{3\mbox{\scriptsize{inj}}(y)}(y,E))$ is a free abelian
group of rank one or two. This group is called the local fundamental
group of $y$ and denoted by $\pi_{1,y}$. If the rank of $\pi_{1,y}$ is
one, $B_{\mbox{\scriptsize{inj}}(y)}(y,E)$ is a ball with only
one pair of points identitied on the boundary, i.e., there is
a unique shortest closed geodesic goes through $y$ representing
the generator of the local fundamental group. If the rank of 
$\pi_{1,y}$ is two, then $B_{3\mbox{\scriptsize{inj}}(y)}(y,E)$ contains
the whole cross section which is an essential embedded torus in $E$.

Case 2). $E$ is a boundle over $S^1$ whose fiber is $\mathbb{E}^2$. If the 
boundle is a product bundle, then injectivity radius is a constant 
function on $E$, local fundamental group $\pi_{1,y}$ has rank one and
there is a unique shortest closed geodesic goes through $y$ representing
the generator of the local fundamental group. Otherwise $E$ is a quotient 
of $\mathbb{E}^3=\mathbb{E}^1 \times \mathbb{E}^2$ by an infinite cyclic group 
generated by $g$ which is equal to a translation by $d$ in $\mathbb{E}^1$
times a rotation with angle $\theta \neq 0$ on $\mathbb{E}^2$. In this
case, the isometry group of $E$ acts transitively on the equi-distant
tori from the unique closed geodesic $\gamma$ whose length is $d$. 
Take a point $y$ on an equi-distant torus with
distance $r \geqslant 60d$. We can find an element $g^i$ ($1\leqslant i \leqslant 20$) 
such that the absolute value of its rotation angle is less than $\pi /10$, 
then we have that
$$
\mbox{inj}(y) \leqslant \frac{1}{2} \left( \frac{\pi r}{10} + 20 d \right) \leqslant 
\frac{r}{3}.
$$
This means that if $y\notin \cup _{y' \in \gamma}
B_{60 d}(y',E)$, 
then $B_{3\mbox{\scriptsize{inj}}(y)}(y,E) \subset  (\mathbb{E}^1 
\times (\mathbb{E}^2 - \{ 0 \})) / \langle g \rangle$. As same as in
Case i) the local fundamental group $\pi_{1,y}$ of $y$ has rank
one or two depends on we have only one or more pairs of points identitied 
on the boundary of $B_{\mbox{\scriptsize{inj}}(y)}(y,E)$. If the
rank is one there is a unique shortest broken closed geodesic 
goes through $y$ representing the generator. If the rank is two, 
then $B_{3\mbox{\scriptsize{inj}}(y)}(y,E)$ contains the whole
equi-distant torus.

Case 3). $E$ is a boundle over $S^1$ whose fiber is an infinite disk
with a cone point.Take a point $y$ on an equi-distant torus with
distance $r$ at least $60$ times of the length of the singular
locus. The same argument as in the Case 2) yields the same conculsion.

For hyperbolic cone manifold $C$, we can argue similarly.
Let $N= \cup _{x\in \Sigma} B_{60\mbox{\scriptsize{inj}}(x)}(x,C)$,
this is a disjoint union of solid torus with a singular core, because 
that $B_{120\mbox{\scriptsize{inj}}(x)}(x,C) \cap \Sigma$ has only
one sigular core for any $x \in \Sigma$. Now let $D=C-N$.
For any point $x\in D$, $B_{3\mbox{\scriptsize{inj}}(x)}(x,C)\cap 
\Sigma =\varnothing$. Define the image of $\pi_1
(B_{\mbox{\scriptsize{inj}}(x)}(x,C)) \subset
\pi_1(B_{3\mbox{\scriptsize{inj}}(x)}(x,C))$ the local fundamental
group $\pi_{1,x}$ of $x$, which is a free abelian group of rank one
or two depends on there is one or more pairs of points identified on
the boundary of $B_{\mbox{\scriptsize{inj}}(x)}(x,C)$. If the rank
is one, there is a unique shortest broken closed geodesic 
goes through $x$ representing the generator of $\pi_{1,x}$.
If the rank is two, the holonomy image of $\pi_{1,x}$ is abelian,
this implies that the holonomy image of $\pi_{1,x}$ is generated 
by one hyperbolic element or two loxodromic elements with a common 
axis because we have no parabolic elements if $C$ is compact. Either
cases, $B_{3\mbox{\scriptsize{inj}}(x)}(x,C)$ contains the whole
equi-distant torus from the singular circle or the closed geodesic
with length $\leqslant \delta$ and once a point on a such equi-distant 
torus has rank two local
fundament group, all local fundamental group of other points on 
the same equi-distant torus have rank two as well.

Consider all a tube which is in form of $T_{\gamma, r} = \cup_{y \in 
\gamma} B_r(y, C)$ and $\gamma$ is a component in $\Sigma$ or a 
closed geodesic with length at most $\delta$ with the property that
any point on the equi-distant torus has rank two local fundamental 
group. All these $\gamma$ has only finitely many choices. For a
fixed $\gamma_1$, if there exists a such $T_{\gamma_1, r}$, we can
define a maximal one. Let $\{ T_i \}$ be all these maximal tubes, we
can see that they are pairs-wisely disjoint. Otherwise $T_i \cap
T_j \neq \varnothing$, then $\partial T_i$ and $\partial T_j$ are parallel.
This shows that $X- \Sigma$ is either a solid torus or a product of
$\mathbb{E}^1$ and a torus, contradict to theorem 4.1.

Denote $D' = D - \cup_i T_i$ which is homeomorphic to $C$ minus those
tubes whose core are not a singular componenmt. Now for any point
$x\in D'$, the local fundamental group $\pi _{1,x}$ is
isomorphic to $\mathbb{Z}$. The generator $\alpha _x$ of $\pi _{1,x}$ is
realized by the unique shortest broken geodesic through $x$, If we fix 
a point $x_0$, take any closed curve $c$ with base point $x_0$ and move the 
unique shortest broken geodesic through $x$ along $c$, we will end up 
with the same broken geodesic. This shows that 
$\pi _{1,x_0}$ is a normal subgroup of $\pi _1(D')$. We know that $D'$ 
is irreducible and we also know that $\partial D' \neq \varnothing $, 
so $D'$ is Haken. These two facts together imply that $D'$ is Seifert
fibered and $\pi_{1,x_0}$ is carried by the fiber.

If $T_i$ is a maximal tube with core is a closed geodesic, then the 
1-dimensional fiber on the boundary of $T_i$ is not a meridian of 
$T_i$, since the fiber represents a nontrivial element in the local 
fundamental group and it has nontrivial holonomy image. So we can extend
the foliation on $D'$ to a Seifert fiberation on $X-\Sigma $. Again
this is in contradiction with Theorem 4.1, so we complete the proof 
of the theorem. \hspace*{\fill}$\square $

{\sc Remark.} In the proof of the theorem, we first rule out the
case that sequences of rescaled hyperbolic cone manifolds converge
to a compact Euclidean cone manifold, then argue that there are only 
three possible noncompact Euclidean manifolds can occur in the limit, 
which are a product of $\mathbb{E}^1$ with a torus, a bundle over $S^1$ 
whose fiber is $\mathbb{E}^2$ or an infinite disk with a cone, all by the
condition $\omega < 2 \pi /3$. If these things can be done then the
rest of the proof goes through. For example, we can also prove the
following theorem.

{\sc Theorem 4.3.} {\em Let $\omega < \pi$, $(X, \Sigma)$ a 
combinatorial type of a hyperbolic cone structure, $X$ is 
not a spherical manifold and there is no embedded sphere 
$S^2 \subset X$ which intersects $\Sigma$ transversely, 
so that $S^2 \cap \Sigma$ is the set of three points. Then
there is a constant $\delta >0$, for any 3-dimensional hyperbolic 
cone structure $C$ on $(X, \Sigma)$ with all cone angles 
less than $\omega $, we have $C_{ \mbox{\scriptsize{thick, }}\delta }
\neq \varnothing $.}

{\sc Proof.} If one sequence of rescaled hyperbolic cone structures 
on $(X, \Sigma)$ converges to a compact Euclidean cone structure,
$(X, \Sigma)$ is a combinatorial type of a compact Euclidean
cone manifold. We know that $X$ is either a spherial manifold or
$X - \Sigma$ is Seifert fibered. Both are impossible. For the noncompact
Euclidean cone manifold, which may occur in the limit, we have one 
more case to consider, that is a product of $\mathbb{E}^1$ with a double of
an acute angled Euclidean triangle. If it can happens, by Theorem 1.3, 
we can find an embedded sphere $S^2 \subset X$ which intersects 
$\Sigma$ transversely, so that $S^2 \cap \Sigma$ is the set of 
three points. That is also ruled out by our hypothesis. The rest of the 
proof is the same as the proof of Theorem 4.2.
\hspace*{\fill}$\square$

\vspace{0.1in}

{\bf \S 5 Compactness of Hyperbolic Cone Structure}

In this section, we will prove a compactness theorem for hyperbolic cone
structures. It is the key step toward Theorem A.

For a fixed combinatorial type $(X,\Sigma )$, let $K$ be a
triangulation of $(X,\Sigma )$ such that the simplicial neighborhood 
of $\Sigma $ in $K$ is a regular neighborhood of $\Sigma$. Note that 
this is always true, if we replace $K$ by its second barycentric derived
subdivision. Then the volume of hyperbolic cone manifolds with the same
combinatorial type $(X, \Sigma )$ has a uniform bound.

{\sc Theorem 5.1.} ([Proposition 4.1, {\bf SOK}]) {\em For all 
3-dimensional hyperbolic cone manifolds with the same combinatorial 
type $(X,\Sigma )$, the volumes have an upper bound $\nu n_3$, 
where $\nu = -3 \int_0^{\pi /3} \log | \sin 2u | du$ is the
maximal volum of a hyperbolic simples (see \em[{\bf M}]\em) and
$n_3$ is the number of 3-simplexes in $K$.}

{\sc Corollary 5.2.} {\em Let $C$ be a compact 3-dimensional 
hyperbolic cone manifold, $\delta >0$ be any positive number. Then 
the diameter of any component of $C_{\mbox{\scriptsize{thick,}}
\delta }$ has a uniform upper bound which only depends on the 
combinatorial type of $C$, $\delta $ and the lower bound of
cone angles $\omega$.}

{\sc Proof.} Let $V$ be the upper bound of the volume, $\omega$ 
the lower bound of the cone angles and $C_1$ a component of $C_{\mbox{\scriptsize{thick, }}\delta }$. For each point $x \in 
C_1 \cap \Sigma$, $B_{\delta /4}(x)$ is standard. We put disjoint
standard $\delta /4$ balls centered in $C_1 \cap \Sigma$ as many as 
we can and denote the cenetrs of the balls by $\{ x_i \}$. Note 
that $N_{\delta /2} (\Sigma) = \{ x \in C_1 |
d(x, C_1 \cap \Sigma) \leqslant \delta /2 \} \subset \cup_i B_{\delta}
(x_i)$. So $\{ x_i \}$ is a $\delta$-net of $N_{\delta /2} ( \Sigma)$.
Then for those points $x \in C_1 - N_{\delta /2}$, balls centered at
these points with radius $\frac{\delta}{2} \sin \frac{\omega}{2}$
are standard. Use these balls to packing $C_1 - N_{\delta /2 }( \Sigma)$,
and the centers of these balls together with $x_i$'s form a 
$\delta $-net of $C_1$. Total number of points in this 
$\delta $-net is bounded, so diameter of $C_1$ is bounded.
\hspace*{\fill}$\square $

Another theorem associated with the volume is the following:

{\sc Theorem 5.3.} {\em The volume of hyperbolic cone manifolds is 
lower semicontinuous, i.e., if $(C_n,x_n)$ is a convergent sequence 
and the limit $(C,x)$ is still a cone manifold, then}
$$
\mbox{vol}(C,x)\leqslant \lim \inf \mbox{ vol}(C_n,x_n). 
$$

{\sc Proof.} We can suppose that $\lim $ vol$(C_n,x_n)$ exists, say it is
equal to $v$. We only need to show that vol$(B_R(x,C))\leqslant v$ for all $R$. 
$B_R(x,C)=\lim B_R(x_n,C_n)$ and the $B_R(x,C)$ and the $B_R(x_n,C_n)$ are
compact. So vol$(B_R(x))=\lim $ vol$(B_R(x_n))$, which is less than $v$.
\hspace*{\fill}$\square $

{\sc Corollary 5.4.} {\em Let $(C_n,x_n)$ be a convergent sequence of
compact hyperbolic cone manifolds with fixed combinatorial type, if the
limit is still a hyperbolic cone manifold, then $\lim (C_n,x_n)$ has finite
volume.}

{\sc Proof.} Applying Theorem 5.1, we get that $\lim \inf $ vol$
(C_n,x_n)<\infty $. Then the conclusion follows from Theorem 5.3. 
\hspace*{\fill}$\square $

The main theorem in this section is the following:

{\sc Theorem 5.5.} {\em For any $0<\varepsilon <\omega <2\pi /3$, the 
space of all 3-dimensional hyperbolic cone structures with the same 
combinatorial type $(X,\Sigma )$ and all cone angles between 
$\varepsilon $ and $\omega $ is compact under Hausdorff distance.}

{\sc Proof.} Let $C_n$ be a sequence of compact 3-dimensional hyperbolic
cone manifolds with the same combinatorial type $(X,\Sigma )$ and all 
cone angles between $\varepsilon $ and $\omega $. By Theorem 4.2, there 
is a $\delta >0$, so that $C_{n,\mbox{\scriptsize{thick,}}\delta }
\neq \varnothing$. Pick $x_n \in C_{n,\mbox{\scriptsize{thick,}}\delta }$ 
and consider the sequence $(C_n,x_n)$, by Theorem 1.2, it has a convergent
subsequence and the limit $(C,x)$ of this subsequence is still a hyperbolic
cone manifold. If $C$ is compact, it will have the same combinatorial type
as $C_n$ by Theorem 1.3. This is the conclusion of our theorem.

Suppose that $C$ is noncompact and we will show that it leads 
to a contradiction. Corollary 5.4 says that $C$ has finite 
volume. Without loss of generality, we also can assume that 
$(C_n,x_n)$ is convergent itself. The number of circle components 
in $\Sigma _C$ is less than the number of components in $\Sigma $. 
We can find a number $M_0>0$, so that all the circle components in 
$\Sigma _C$ are contained in $B_{M_0}(x,C)$. The rest of the proof 
will be divided into two cases.

Case 1). Suppose that $[C-B_m(x,C)]\cap \Sigma _C\neq \varnothing $ for all 
$m>M_0$. Let $y_m\in [C-B_m(x,C)]\cap \Sigma _C$, then $\lim _{m\rightarrow
\infty }$inj$(y_m)=0$, as $C$ has finite volume. For $m$ large enough and
applying Theorem 1.4, we know that $C$ looks like a noncompact Euclidean
cone manifold $E$ near $y_m$. Since all cone angles are less than $\omega $, 
$E$ must be a bundle over $S^1$ whose fiber is an infinite disk with a cone.
In the Euclidean case, the singular component is a circle, so the singular
component which contains $y_m$ is a circle too and it contradicts the
hypothesis that all circle components of $\Sigma _C$ are contained in 
$B_{M_0}(x,C)$.

Case 2). $[C-B_m(x,C)]\cap \Sigma _C=\varnothing $ for some $m>M_0$. This
really means that all singular component are circles. Then, far 
away from $x$, $C$ looks like a noncompact Euclidean manifolds. By 
the same reason as in the proof of Theorem 4.2, $E$ can not be a 
twisted line bundle over a Klein bottle. So we can do thin-thick 
decomposition for $C$, and all ends of $C$ are cusps.

Let $(N, \Sigma_C)$ be a compact core of $(X_C,\Sigma _C)$. By 
Theorem 1.3, we can 
embed $N$ back into $C_n$ for $n$ large enough. Let $g_n$ be the 
embeddings and $N_n$ the image $g_n(N)$. $\partial N_n$ is a bunch of 
tori, and each torus either bounds a solid torus in $X-\Sigma $ or 
is parallel to a singular circle. For each boundary torus $T$ of $N$, 
we have a sequence $\alpha_n^T \in {\bf H}_1(T)$, where $\alpha_n^T$ 
is represented by the preimage of the essential circle $\beta _n$ in 
$g_n(T)$ which is null homotopic or a meridian. $\alpha_n^T$ must tend 
to infinity in ${\bf H}_1(T)$. Otherwise, $\alpha _n^T$ has a subsequence, 
we still call it $\alpha _n^T$, which converges to $\alpha $. The 
holonomy image of $\beta _n$ will converge to the holonomy image of 
$\alpha $, which is a parabolic element. The holonomy of $\beta _n$ 
is either a trivial element or an elliptic element with rotation 
angle between $\varepsilon $ and $\omega $. Such a sequence cannot
converge to a parabolic element.

Now we can apply Thurston's hyperbolic Dehn surgery theorem to this
situation. For $n$ large enough, the Dehn filling along all these 
$\alpha_n^T$'s will give us a hyperbolic manifold $\hat N_n$. Note 
that all these results of Dehn filling are the same, they are all 
homeomorphic to $X$. On the other hand, the volume of a hyperbolic 
manifold is an invariant and vol$(N)>$vol$(X)=$vol$(\hat N_n)
\rightarrow $vol$(N)$ (see Theorem 6.5.6, [{\bf T1}] and Theorem 1, 
[{\bf NZ}]). This is impossible.

So $C$ must be compact, and we finish the proof of the theorem. 
\hspace*{\fill}$\square $

The same as the Remark after Theroem 4.2 and the proof of Theorem
4.3, we can prove the following theorem.

{\sc Theorem 5.6.} {\em Let $0< \varepsilon < \omega < \pi$, 
$(X, \Sigma)$ a combinatorial type of a hyperbolic cone structure, 
$X$ is not a spherical manifold and there is no embedded sphere 
$S^2 \subset X$ which intersects $\Sigma$ transversely, 
so that $S^2 \cap \Sigma$ is the set of three points. Then the 
space of all 3-dimensional hyperbolic cone structures on $(X,\Sigma )$ 
with all cone angles between $\varepsilon $ and $\omega $ is 
compact under Hausdorff distance.}

\vspace{0.1in}

{\bf \S 6 Deformation Theory of Hyperbolic Cone Structures}

We will discuss the deformation of hyperbolic cone structures in this
section. Most of the theorems in this section are known to Thurston.

{\sc Theorem 6.1.} {\em Let $C$ be a hyperbolic cone manifold. The
holonomy represetation $\rho: \pi_1(X_C-\Sigma_C) \rightarrow
PSL_2(\mathbb{C})$ can be lifted to a representation in $SL_2(
\mathbb{C})$.}

{\sc Proof.} Our proof goes the same fashion as the proof of
[Proposition 3.1.1, {\bf CS}]. Let $\rho: \pi_1(X_C-\Sigma_C) 
\rightarrow PSL_2(\mathbb{C})$ be the representation, then $\pi_1
(X_C-\Sigma_C)$ is isomorphic to the subgroup $\Gamma = \{(g, 
\rho(g))\}$ in $\pi_1(X_C-\Sigma_C) \times PSL_2(
\mathbb{C})$. $\Gamma$ acts on $(\widetilde{X_C-\Sigma_C}) \times
PSL_2(\mathbb{C})$ freely and proper discontinuousely, the quotient 
$Q$ is a principle $PSL_2(\mathbb{C})$ boundle over $X_C-\Sigma_C$.
$PSL_2(\mathbb{C})$ can be identified with the principal bundle of
orthonormal tangent frames to $\mathbb{H}^3$. Use the fact that
$X_C-\Sigma_C$ has trivial tangent bundle, it is not difficult 
to see that $Q$ is homeomorphic to $(X_C-\Sigma_C) \times \mathbb{H}^3
\times SO(3)$. 

On the other hand, let $\tilde{\Gamma}$ is the 
preimage of $\Gamma$ under the homomorphism $\pi_1(X_C-
\Sigma_C) \times SL_2(\mathbb{C}) \rightarrow \pi_1(X_C-\Sigma_C) 
\times PSL_2(\mathbb{C})$, then $Q = (\widetilde{X_C-\Sigma_C}) \times
SL_2(\mathbb{C})/\tilde{\Gamma}$. $SL_2(\mathbb{C})$ is simply connected,
so $\tilde{\Gamma} = \Gamma \times \mathbb{Z}_2$. This shows that
$\rho$ can be lifted to a representation in $SL_2(\mathbb{C})$.
\hspace*{\fill}$\square$

Later on, the lift $\rho: \pi_1(X_C-\Sigma_C) \rightarrow
SL_2(\mathbb{C})$ in Theorem 6.1 will be also called holonomy
representation. 

Let $(X, \Sigma)$ be a combinatorial type of a compact 
hyperbolic cone structure, $\Sigma$ has $m$ components 
$S_1, \ldots, S_m$ and let $\mu_i$ be the meridian of $S_i$.

{\sc Theorem 6.2.} {\em Let $C_0$ be a hyperbolic cone manifold 
structure on $(X,\Sigma )$ or a complete hyperbolic structure on 
$X-\Sigma $ and $\rho_0$ its holonomy representation. If $\rho$ 
is sufficiently close to $\rho_0$ such that all $\rho (\mu_i)$
are elliptic. Then there is a hyperbolic cone manifold structure 
$C_1$ on $(X,\Sigma )$ its holonomy representation is $\rho$.
Furthermore, if $C_0$ is a hyperbolic cone structure, then $C_1$
can be chosen close to $C_0$ and such nearby cone structure
is determined uniquely by the conjugacy class of $\rho$.}

{\sc Proof.} Let ${\mathcal N}(\Sigma )$ be a regular neighborhood 
of $\Sigma $. Then $U=X-{\mathcal N}(\Sigma )$ is compact. 
By [Proposition 5.1, {\bf T1}], there is a hyperbolic 
structure $C_1^{\prime }$ on $U$ which is close to $C_0|_U$ with 
holonomy $\rho$. Since each $\rho(\mu _i)$ are elliptic, 
we have a unique way to extend $C_1^{\prime }$ to a hyperbolic 
cone structure $C_1$ on $(X,\Sigma )$. This $C_1$ is determined 
uniquely by the conjugacy class of $\rho$. \hspace*{\fill}$\square$

The representation space plays an essential role in the 
deformation theory for hyperbolic cone structures. We denote 
by $R(X,\Sigma)$ the representation space $R(X, \Sigma) = 
\{\rho :\pi _1(X-\Sigma ) \rightarrow SL_2(\mathbb{C})\}$, which 
is a complex affine algebraic set. When we consider different hyperbolic 
cone structures on a fixed combinatorial type $(X, \Sigma)$, 
define marked hyperbolic cone structure is convenient. 
A marked hyperbolic cone structure on $(X,\Sigma)$ is an 
equivalent class of hyperbolic cone structures on $(X,\Sigma )$ 
with equivalent relation that there is an isometry which is 
homotopic rel$\Sigma$ to the identity map and it is similar 
to marked Riemann surface in Teichm\"uller theory. Theorem 
6.2 says that near a marked hyperbolic cone structure, the 
space of marked hyperbolic cone structures is parametrized 
by conjugacy classes of their holonomy representations.
A deformation of conjugacy class of the holonomy representation of 
a cone structure in $R(X, \Sigma)/SL_2(\mathbb{C})$ is called 
an algebraic deformation of the cone structure. locally a 
geometric deformation is equivalent to an algebraic deformation, 
and increasing cone angles is the same as increasing the 
angle of elliptic elements $\rho_0(\mu_i)$.

The basic fact about the algebraic deformation is the following
theorem.

{\sc Theorem 6.3.} {\em Let $R_0$ be an irreducible component of
$R(X,\Sigma)$ which contains $\rho_0$ a holonomy representation of
a compact hyperbolic cone structure with combinatrial type 
$(X, \Sigma)$ or the complete hyperbolic structure on $X-\Sigma$, 
then 
$$
\dim_{\mathbb{C}} R_0 \geqslant m + 3,
$$
where $m$ is the number of circle components of $\Sigma$.}

{\sc Proof.} $\rho_0$ is irreducible and let $N$ be a compact 
core of $X-C$, then for each torus component $T$ 
of $\partial N$, $\rho_0( \pi_1 (T) ) \not\subset \{
\pm 1 \}$. [Propersition 3.2.1, {\bf CS}] says that 
$\dim_{\mathbb{C}} R_0 \geqslant m + 3$. \hspace*{\fill}$\square$

For a compact marked hyperbolic cone structure, the holonomy 
repersentation is always irreducible and unique upto conjugacy. 
So the space of characters will be more interesing, it is because 
the conjugacy classes of irreducible representations are 
one-to-one corresponding to their characters. Let $R_0$ be an 
irreducible component of $R(X,\Sigma)$ which contains $\rho_0$ 
a holonomy representation of
a compact hyperbolic cone structure with combinatrial type 
$(X, \Sigma)$ or the complete hyperbolic structure on $X-\Sigma$,
${\mathcal C}_0$ the space of all characters of representations in $R_0$.

{\sc Proposition 6.4.} {\em ${\mathcal C}_0$ is an affine variety 
and $\dim_{\mathbb{C}}{\mathcal C}_0=\dim_{{\mathbb{C}}}R_0 - 3 \geqslant m$.}

See [Section 1, Proposition 3.2.1, {\bf CS}] for a detailed proof.

Let $f:{\mathcal C}_0 \rightarrow \mathbb{C}^m$ be the trace map defined by
$$
f([\rho ])=(\mbox{tr}(\rho (\mu _1)),\ldots ,\mbox{tr}(\rho (\mu _m))). 
$$
$f$ is a regular map and $f([\rho _0])=(\epsilon _12\cos \theta _1,
\ldots, \epsilon _m2\cos \theta _m)$, where $\epsilon _i=\pm 1$, 
since each $\rho_0(\mu _i)$ is elliptic with a rotation angle 
$\theta _i$. 

Suppose that $\rho_0$ is a holonomy representation of 
the complete hyperbolic structure on $X-\Sigma$, Mostow's rigidity
theorem says that $f^{-1}(f([\rho _0]))$ is a single point in 
${\mathcal C}_0$, which implies that $\dim_{\mathbb{C}}{\mathcal C}_0 = m$ and 
the trace map $f$ is onto near $[\rho_0]$.
For the case that $\rho_0$ is a holonomy representation of a compact 
hyperbolic cone structure with all cone angles less than $2 \pi /3$, 
the same conclusion holds. 

{\sc Proposition 6.5.} {\em Suppose $\theta _i<2\pi /3$ for all $1 \leqslant 
i \leqslant m$, then the irreducible component of $f^{-1}(f([\rho _0]))$ in 
${\mathcal C}_0$ which contains $[\rho _0]$ is the single point $[\rho _0]$.}

{\sc Proof.} Let $L$ be an irreducible component of $f^{-1}(f([\rho _0]))$
which contains $[\rho _0]$ and $S$ be the subset of $L$ which consists of
all characters which can be realized by a holonomy representation of
hyperbolic cone structures. Theorem 6.2 says that $S$ is open in $L$. 
On the other hand, the convergence under the Hausdorff distance implies 
the convergence of the holonomy representations, so Theorem 5.5 says that
$S$ is also compact in $L$. In particular, $S=L$ is a compact irreducible 
affine variety, and thus it must be the single point $[\rho _0]$.
\hspace*{\fill}$\square$

The proposition implies that $\dim_{\mathbb{C}} {\mathcal C}_0 = m$ and
$f$ is onto near $[\rho_0]$. Now we restrict our attention to the 
case that all $\theta _i < 2\pi /3$. For any given direction 
$\vec a=(a_1, \ldots ,a_m)$, we define a ray $l(t)=f([\rho_0])+
t\vec a$ in ${\mathbb{C}}^m$.

{\sc Theorem 6.6.} {\em There is a continuous lift of $l(t)$ 
$$
\tau :[0,T)\rightarrow {\mathcal C}_0, 0< T < \infty , 
$$
such that $\tau (0)=[\rho _0]$.}

{\sc Proof.} Let $L \subset \mathbb{C}^m$ be the complex line which 
contains $l(t)$. Propositions 6.4 and 6.5 then implies that the 
irreducible component $E$ of $f^{-1}(L)$ which contains $[\rho _0]$ 
must be an affine curve in ${\mathcal C}_0$ and $f |_{E}$ is a nonconstant
map. Any a such $f|_E : E\rightarrow L$ is open, so we can lift $l(t)$ to 
$\tau :[0,T) \rightarrow {\mathcal C}_0$. \hspace*{\fill}$\square $

{\sc Remark.} Proposition 6.5 and Theorem 6.6 also hold, if we release
the restriction $\theta_i < 2 \pi /3$ to $\theta_i < \pi$ and
put additional condition on $(X, \Sigma)$ as Theorems 4.3 or 5.6.

Theorem 6.6 tell us that local algebraic deformations 
of cone structure always exist, so do local geometric deformations 
by Theorem 6.2. The complete hyperbolic structure can be regarded 
as a cone structure on $(X, \Sigma)$ with ``zero cone angle''. This 
theorem is also true when we replace the $\rho_0$ by the representation 
of the complete hyperbolic structure on $X-\Sigma $. Instead of using 
Proposition 6.5, we can use Mostow's rigidity theorem in the proof. 
Everything works fine. However, if we still want to get a geometric 
deformation of hyperbolic cone structure, we need assume that 
$a_i \epsilon_i \leqslant 0$ for all $i$. This will guarantee that all 
holonomy image of $\mu _i$ are elliptic, therefore we can get 
the inverse of Theorem 4.1. We leave the details of the proof to the
reader.

{\sc Theorem 6.7.} {\em Suppose that $\Sigma$ is a hyperbolic link in 
a closed mainfold $X$. There is an $\varepsilon >0$, for any 
$\theta = ( \theta_1, \ldots, \theta_m) \in \mathbb{R}^m$ and 
$0< \theta_i < \varepsilon$, we have a hyperbolic cone structure 
$C_{\theta}$ on $(X,\Sigma )$ such that the cone angle of $S_i$ 
is $\theta_i$.}

To finish the section, we prove the following lemma, which says
that, just like the complete hyperbolic structures, the 
orbifold structures are also rigid.

{\sc Lemma 6.8.} {\em Let $O_1$ and $O_2$ be two hyperbolic orbifold
structures on an orbifold $O$. then they are isomorphic. In fact the
isometry can be chosen to be homotopic to the identity map rel$\Sigma$.}

{\sc Proof.} Any hyperbolic orbifold is very good, so there is 
a finite regular branch cover $\tilde O$ of $O$ which is a manifold. 
Denote by $G$ the covering transformation group. Pulling back the 
structures $O_1$ and $O_2$ to $\tilde O$, we get two complete hyperbolic 
structures $\tilde O_1$ and $\tilde O_2$ on $\tilde O$ such that 
elements of $G$ are isometries of both structures $\tilde O_i$. 
Mostow's rigidity theorem says that we can homotopy the identity 
on $\tilde O$ to an isometry $\tilde h : \tilde O_1 \rightarrow 
\tilde O_2$. For any element $g \in G$, $g$ and $\tilde h^{-1} g
\tilde h$ are two homotopic isometries on the structure $\tilde O_1$. 
This implies that $g = \tilde h^{-1} g \tilde h$,
i.e., $\tilde h$ is $G$-invariant. So $\tilde h$ induces an isometry 
$h: O _1 \rightarrow  O_2$.\hspace*{\fill}$\square $

\vspace{0.1in}

{\bf \S 7 Moduli Space of Hyperbolic Cone Structures}

In this section, we will prove Theorem A and Corollary B mentioned 
in the begining of the paper.

Denote by $\rho_*$ the representation of the complete hyperbolic 
structure on $X-\Sigma $, by $R_*$ the irreducible component of
$R(X, \Sigma)$ which contains $\rho_*$ and ${\mathcal C}_*$ the space
of characters of representations in $R_*$. 
Let $C$ be another hyperbolic cone structure on $(X,\Sigma)$, 
it has a holonomy $\rho \in R(X,\Sigma )$. Let $R_\rho$ be the
irreducible component of $\rho$ in $R(X,\Sigma )$ and 
${\mathcal C}_\rho $ the corresponding character space. If all cone 
angles are less than $2 \pi /3$, we have the following

{\sc Lemma 7.1.} {\em $R_\rho = R_*$ and ${\mathcal C}_\rho ={\mathcal C}_*$.}

{\sc Proof.} It suffices to prove that $R_\rho = R_*$. 
By Theorem 6.7, we can choose a hyperbolic orbifold structure 
$O$ on $(X, \Sigma)$ with all cone angles are $2 \pi /n$ and 
its holonomy representation lies in $R_*$. Let
$$
l(t)=(1-t)(\epsilon _12\cos \theta _1,\ldots ,\epsilon _m2\cos \theta _m)+t(\epsilon _12\cos 2\pi /n,\ldots ,\epsilon_m2\cos 2\pi /n), 
$$
be a line in $\mathbb{C}^m$ and $I \subset [0, 1]$ be the largest interval 
on which there is a lift $\tau: I \rightarrow {\mathcal C}_\rho$ of $l(t)$ 
such that $\tau (0) = [\rho]$ and each $\tau (t)$ is the charactor of 
a holonomy representation of a hyperbolic cone structure on $(X, \Sigma)$. 
Theorems 6.2 and 6.6 imply that $I$ is open in $[0,1]$, and Theorem 5.5 
says that $I$ is also compact. Therefore, it is whole interval 
$[0,1]$. 

$\tau (1)$ is the charactor of a holonomy representation $\rho'$ of 
a hyperbolic cone structure on $(X, \Sigma)$. $l(1) = (\epsilon _1
2\cos 2\pi /n, \ldots, \epsilon_m2\cos 2\pi /n)$ implies that all 
cone angles of this structure are $2 \pi /n$, so it is a hyperbolic 
orbifold. Denote this structure by $O'$ and note that the holonomy 
representation $\rho' \in R_{\rho}$. Apply
Lemma 6.8, we know that $O$ and $O'$ are isometric by an isometry
homotopic to the identity map rel$\Sigma$, this means that they have
the same holonomy representation, i.e., $\rho' \in R_{\rho} \cap
R_*$. The irreducible component $R_0$ which contains $\rho'$ is
contained in $R_{\rho} \cap R_*$, and all these three affine vatieties
$R_0, R_{\rho}$ and $R_*$ have dimension $m$, so we have 
$R_0 = R_{\rho} = R_*$. \hspace*{\fill}$\square $

Now we are in the position to prove that the trace map $f: {\mathcal C}_*
\rightarrow \mathbb{C}^m$ is a local homeomorphism near any
character of a holonomy representation of a hyperbolic cone
structure.

{\sc Theorem 7.2.} {\em Let ${\mathcal U}$ be the subset of ${\mathcal C}_*$
consists of all characters of holonomy representations of 
hyperbolic cone structures on $(X, \Sigma)$ with all cone angles
less than $2\pi/3$. Then the trace map $f: {\mathcal U} \rightarrow 
\mathbb{C}^m$ is a homeomorphism onto its image.}

{\sc Proof.} Let ${\mathcal W} = \{(\epsilon _12\cos \theta _1,
\ldots ,\epsilon _m2\cos \theta _m) | 0< \theta_i < 2 \pi /3 \}
\subset \mathbb{C}^m$. By the proof of Lemma 7.1, we can see that 
the trace map $f({\mathcal U}) = {\mathcal W}$ is onto.

To show that $f: {\mathcal U} \rightarrow {\mathcal W}$ is a
homeomorphism. Let $B_0 \subset {\mathcal W}$ be the set
$$
\{p |f^{-1}( p ) \cap {\mathcal U}\mbox{ contains 
at least two different characters}\} 
$$
and $B=\overline{B}_0$, the closure of $B_0$ in ${\mathcal W}$. We want to 
show that $B$ is open in ${\mathcal W}$. If it is true, $B$ is either empty 
or the whole ${\mathcal W}$. If $B={\mathcal W}$, all points of the form 
$\{(\epsilon _12\cos 2\pi /n_1,\ldots ,\epsilon _m 2\cos 2\pi /n_m)\}
\subset B-B_0$, since they are corresponding to orbifold structures. So 
any character of holonomy representation of a orbifold structure is either
a singular point of ${\mathcal C}_0$ or a singular point of $f:{\mathcal C}_0 
\rightarrow \mathbb{C}^m$. This contradicts to the facts that the image 
of singular set of ${\mathcal C}_0$ and singular set of $f$ lies in a 
lower dimensional algebraic set and $\{(\epsilon_1 2\cos 2\pi /n_1,
\ldots,\epsilon _m 2\cos 2\pi /n_m)\}$ is a Zariski dense set in 
$\mathbb{C}^m$. So $B$ is empty, i.e., $f$ is injective. For any compact 
set $K \subset {\mathcal W}$, the space of hyperbolic cone structures whose
the character of the holonomy representation is contained in 
$f^{-1}(K)$ is compact by Theorem 5.5, so $f^{-1}(K)$ is compact and
this shows that $f^{-1}: {\mathcal W} \rightarrow {\mathcal U}$, the inverse 
of $f$ is also continuous.

Now we prove that $B$ is open. For any $p \in B_0$, $f^{-1}(p) \cap
{\mathcal U}$ contains two distinct characters $\chi^1$ and $\chi^2$ both
corresponding to hyperbolic cone structures. Theorem 6.2 says that 
there are two small neighborhoods $\chi^i \in U^i$ such that all points 
in $U^i\cap f^{-1}({\mathcal W}) \subset {\mathcal U}$. Without loss of 
generality, we can assume that $U^1 \cap U^2 = \varnothing$. Also we know 
that $f$ is onto near a cone structure by Theorem 6.6. Then $f(U^1)\cap f(U^2)\cap \mathbb{R}^m \subset B_0$ is a neighborhood of $p$. If $p\in 
B-B_0$, there is a sequence $p_n\in B_0$ such that $p_n\rightarrow p$. 
For each $p_n \in B_0$, we have two distinct characters $\chi_n^1$ and 
$\chi_n^2$ they are corresponding to hyperbolic cone structures
$C_n^1$ and $C_n^2$, and $f(\chi_n^i) = p_n$. Without loss of generality, 
we can assume that both $C_n^1$ and $C_n^2$ converge. So two sequences
$\chi_n^i$ and $\chi_n^2$ converge. Let $\chi_n^i \rightarrow \chi^i$, 
$f(\chi^i) = p$ and $p \not\in B_0$ imply that $\chi^1 = \chi^2 =\chi$.
This shows that $\chi$ is singular point of ${\mathcal C}_0$ or a singular
point of $f$ and the local degree of $f>1$. All points of an open 
neighborhood of $p$ except those points 
in a lower dimensional algebraic set are contained in $B_0$. So this 
open set is a neighborhood of $p$ in $B$. This shows that $B$ is open.
\hspace*{\fill}$\square $

The moduli space of marked hyperbolic cone structures with 
combinatorial type $(X,\Sigma)$ and cone angles less than 
$2\pi /3$ will be denoted by ${\mathcal HC}(X,\Sigma)$.

Let $h: {\mathcal HC}(X, \Sigma) \rightarrow {\mathcal U}$ defined by send
each marked hyperbolic cone structure to the character of its holonomy
representation. By the definition, $h$ is onto, Theorem 6.2 tell 
us that $h$ is a local homeomorphism and Theorem 5.5 implies that 
$h$ has the path lifting property. So $h$ is a covering map and then
we get the following corollary, as $\mathcal U$ is simply connected.

{\sc Corollary 7.3.} {\em $fh: {\mathcal HC}(X, \Sigma) \rightarrow
{\mathcal W}$ is a homeomorphism.}

{\sc Remark.} The map
$$
\begin{array}{rcl}
\sigma :(0,2\pi /3)^m & \rightarrow & {\mathcal W} \\ 
(\theta _1,\ldots ,\theta_m) & \longmapsto & (\epsilon _1
2\cos \theta _1,\ldots ,\epsilon _m2\cos \theta _m) 
\end{array}
$$
is a homeomorphism, so $\iota = h^{-1} f^{-1} \sigma :(0,2\pi /3)^m
\rightarrow {\mathcal HC}(X, \Sigma)$ is a natural parametrization of 
${\mathcal HC}(X,\Sigma)$ by $m$ cone angles. In other words, 
Corollary 7.3 is another version of Theorem A.

As a simple application of Theorem A, we have the following
proposition.

{\sc Proposition 7.4.} {\em Let $\Sigma $ be a hyperbolic link 
in a homological sphere $M^3$. Then $k$-fold cyclic cover of 
$M^3$ branched over $\Sigma $ is a hyperbolic manifold provided 
$k\geqslant 4$.}

{\sc Proof.} Let $N$ be the $k$-fold cyclic cover of $M$ branched 
over $\Sigma $. Then $M$ can be viewed as an orbifold with the 
singular locus $\Sigma $ and all cone angles equal to $2\pi /k$. 
Thus $N$ is a manifold cover of $M$. By Theorem 10.4, the orbifold 
$M$ is hyperbolic, so $N$ is a hyperbolic manifold.
\hspace*{\fill}$\square $

As we said before, if we release the restriction $\theta_i < 2 
\pi /3$ to $\theta_i < \pi$ and put additional condition on 
$(X, \Sigma)$ as Theorems 4.3 or 5.6, the arguements in the proofs of
Theorem 7.2 and Corollary 7.3 work fine. So we have

{\sc Theorem 7.5.} {\em Let $\Sigma$ be a hyperbolic link with 
$m$ components in a 3-dimensional manifold $X$. $X$ is 
not a spherical manifold and there is no embedded sphere 
$S^2 \subset X$ which intersects $\Sigma$ transversely, 
so that $S^2 \cap \Sigma$ is the set of three points. Then the 
moduli space of marked hyperbolic cone structures on the pair 
$(X,\Sigma)$ with all cone angles less than $\pi$ is 
an $m$-dimensional open cube, parameterized naturally
by the $m$ cone angles.}

In the rest of section, we will prove Corollary B which is a 
more general than Proposition 7.4.

{\sc Corollary B.} {\em If $M$ is an irreducible, closed, 
atoroidal 3-manifolds, it is not Seifert manifold and
admits a finite group $G$ action. If the order of $G$ is odd,
the $G$-action is effective and not fixed-point-free, then 
the quotient $M/G$ is a geometric orbifold.}

{\sc Remark.} Pull back the singular geometric structure on 
$M/G$, we will get a geometric structure on $M$. 
The theorem says not only that $M$ has a geometric structure, 
but also $G$ consists of isometries.

This is a special case of Thurston's geometrization theroem, we
refer [Problem 3.46, {\bf K}] for the historical remark 
on Thurston's geometrization theorem. 

A lot of concepts in 3-orbifold theory, such as, irreducibility 
and incompressibility, are very similar to the corresponding 
concepts for manifolds. Following Hodgson in [{\bf Ho}], we 
use {\sc sphere} or {\sc torus} to denote a spherical or a 
Euclidean 2-orbifold and {\sc disk} or {\sc ball} for a quotient 
of a disk or a ball under a finite, orientation preserving linear 
action. A 3-orbifold $O$ is irreducible if every {\sc sphere} 
bounds a {\sc ball}. a 2-suborbifold $F$ is incompressible 
if any 1-suborbifold of $F$ which bounds a {\sc disk} in 
$O-F$ also bounds a {\sc disk} in $F$. An orbifold is atoroidal 
if each incompressible {\sc torus} is $\partial $-parallel. For 
our convenience, we always assume that all 3-orbifolds 
contain no bad 2-orbifolds.

$M/G$ is an orbifold and we denote the orbifold structure be
$O_{M/G}$. Let $X_O$ be the underlying space of $O_{M/G}$ and 
$\Sigma_O$ the singular locus. It is easy to see that $\Sigma_O$
is a link in $X_O$ since any finite group of $SO(3)$ with
odd order is a finite cyclic group. 
A simple geometric argument shows that 
$O_{M/G}$ is closed, irreducible and atoroidal under our assumption
that $M$ is closed, irreducible and atoroidal. Denote $X_O -
\Sigma_O$ by $N_O$. Thurston's program start with the following 
proposition.

{\sc Proposition 7.6.} {\em Either $O_{M/G}$ is Seifert fibered 
or $N_O$ admits a complete hyperbolic structure with finite volume.}

For a proof of the proposition, we refer readers to [\S2, {\bf SOK}].

It is easy to handle the case that $O_{M/G}$ is Seifert fibered. 
Applying the method in [{\bf S}] (also see [{\bf JN}]), we will 
see that all of these orbifolds have Seifert geometries. So we 
assume that $O_{M/G}$ is not a Seifert orbifold and this implies 
that $\Sigma _O \subset X_O $ is a hyperbolic link.

If all cone angles in $O_{M/G}$ are less than $2 \pi /3$, Theorem
A implies that we can find a required hyperbolic orbifold structure on 
$O_{M/G}$. So we now can assume that there exists at least one
cone angle is $2 \pi /3$. In this case, we can find a sequence
of hyperbolic cone structures $C_n$ on $(X_O, \Sigma_O)$ such that
all cone angles tends to desired cone angles. Now we will
investigate the limit of $C_n$. For this purpose, we need following
refinements of Theorems 4.2 and 5.5.

{\sc Proposition 7.7.} {\em For the sequence of hyperbolic
cone structures $C_n$, we have the following alternative:

i) there is a sequence $x_i \in C_n$ and the rescaled sequence 
$( (\mbox{inj}(x_n))^{-1}C_n,x_n)$ has a subsequence 
which converges to a compact Euclidean cone manifold, or

ii) there is a constant $\delta >0$, so that $C_{n,
\mbox{\scriptsize{thick,}}\delta }\neq \varnothing $.}

{\sc Proposition 7.8.} {\em If there is a $\delta>0$, such that
$C_{n,\mbox{\scriptsize{thick,}}\delta }\neq \varnothing $.
$C_n$ has a subsequence $C_{n_k}$ which converges to a
hyperbolic cone structure on $(X_O,\Sigma_O )$.}

{\sc Proof of Propositions 7.7 and 7.8.} Adding the first 
alternative in Proposition 7.7 is really a way to say we 
can assume that there is no rescaled sequence has a subsequence 
which converges to a compact Euclidean cone manifold. 

For the noncompact Euclidean cone manifold, which may occur in 
the limit, just as in the proof of Theorem 4.3, we have one more 
case to consider, that is a product of $\mathbb{E}^1$ with a double of 
an equilateral triangle. If it is the case,
we can embed this {\sc torus} back into $O_{M/G}$. It is easy to see
that this {\sc torus} is incompressible and nonseparating in 
$O_{M/G}$, contradict to the fact that $O_{M/G}$ is atoroidal.
The rest of the proofs are the same as the proofs of Theorems 4.2
and 5.5.
\hspace*{\fill}$\square $

Corollary B is followed immediately by these two propositions.
If the alternative i) of Proposition 7.7 holds, we can get a Euclidean
orbifold structure on $O_{M/G}$; Otherwise, we can use Proposition
7.8 to produce a hyperbolic orbifold structure on $O_{M/G}$. Either
way, $M/G$ is a geometric orbifold. We complete the proof of Corollary B.

\vspace{0.2in}

{\bf References}

\begin{description}

\item [{[{\bf B}]}] A. Beardon, {\em The Geometry of Discrete Groups},
Springer-Verlag New York Heidelberg Berlin 1983.

\item [{[{\bf CS}]}] M. Culler and P. Shalen, Varieties of group 
representations and splitting of 3-manifolds, {\em Ann. of 
Math.} {\bf 117}(1983), 109-146.

\item [{[{\bf G}]}] M. Gromov, Groups of polynomial growth and expanding maps, {\em I.H.E.S. Publications Math.} 
{\bf 53}(1981), 53-78.

\item [{[{\bf Ha1}]}] R. S. Hamilton, Three-manifolds with positive
Ricci curvature, {\em J. Differential Geometry} 
{\bf 17}(1982), 255-306.

\item [{[{\bf Ha2}]}] R. S. Hamilton, Four-manifolds with positive 
curvature operator, {\em J. Differential Geometry} 
{\bf 24}(1986), 153-179.

\item [{[{\bf Ho}]}] C. Hodgson, Geometric structures on 3-dimensional 
orbifolds: Notes on Thurs-ton's proof, preprint.

\item [{[{\bf J}]}] T. Janssen, {\em Crystallographic Groups}, North-Holland Elsevier,
Amsterdam, 1973.

\item [{[{\bf JN}]}] M. Jankins and W. Neumann, {\em Lectures on 
Seifert manifolds}, Brandeis lecture notes 2, March 1983.

\item [{[{\bf K}]}] R. Kirby, Problems in Low-dimensional Topology, 35-474 in:
{\em Geometric Topology, VolII}, ed. W. Kazez, AMS/IP, 1997.

\item [{[{\bf M}]}] J. Milnor, Hyperbolic geometry: the first 150 
years, {\em Bull. Amer. Math. Soc.} {\bf 6}(1982), 9-24.

\item [{[{\bf NZ}]}] W. D. Neumann and D. Zagier, Volumes of hyperbolic
three-manifolds, {\em Topology} {\bf 24}(1985), 307-332.

\item [{[{\bf S}]}] P. Scott, The geometry of 3-manifolds, {\em Bull. 
Lond. Math. Soc.} {\bf 15}(1983), 401-487.

\item [{[{\bf SOK}]}] T. Soma, K. Ohshika and S. Kojima, Towards a 
proof of Thurston's geometrization theorem for orbifolds, 
Res. Inst. Math. Sci., Kokyuroku, {\bf 568}(1985), 1-72.

\item [{[{\bf T1}]}] W. Thurston, {\em The Geometry and Topology of 
Three-manifolds}, Princeton Univ., Math. Dept., 1979.

\item [{[{\bf T2}]}] W. Thurston, Three-manifolds with symmetry, 
preprint (1982).

\item [{[{\bf Z}]}] Q. Zhou, Three-dimensional geometric cone
structures, Ph.D. Thesis, UCLA, 1990.

\end{description}

\end{document}